\documentclass[11pt]{article}
\usepackage{authblk}
\author{Ingvar Ziemann\thanks{Correspondence: \texttt{ingvar.ziemann@gmail.com}}}

\usepackage[utf8]{inputenc} % allow utf-8 input
\usepackage[T1]{fontenc} %\usepackage{ae} \usepackage{aecompl}   % use 8-bit T1 fonts

\usepackage{fullpage}

\usepackage{url}            % simple URL typesetting
\usepackage{booktabs}       % professional-quality tables
\usepackage{amsfonts}       % blackboard math symbols
\usepackage{nicefrac}       % compact symbols for 1/2, etc.
\usepackage{microtype}      % microtypography
\usepackage[dvipsnames]{xcolor}

\usepackage{natbib} %use option numbers if arxiv fails
\usepackage{enumitem}
\usepackage{amsthm}
\usepackage{amssymb}
\usepackage{amsmath}
\usepackage{mathrsfs}
\usepackage{thmtools}
\usepackage{hyperref}
\hypersetup{
    pdftoolbar=true,        % show Acrobat’s toolbar?
    pdfmenubar=true,        % show Acrobat’s menu?
    pdffitwindow=false,     % window fit to page when opened
    pdfstartview={FitH},    % fits the width of the page to the window
    pdfnewwindow=true,      % links in new PDF window
    colorlinks=true,       % false: boxed links; true: colored links
    linkcolor=blue,          % color of internal links (change box color with linkbordercolor)
    citecolor=ForestGreen,        % color of links to bibliography
    filecolor=magenta,      % color of file links
    urlcolor=red,           % color of external links
    breaklinks=false,
}

\usepackage{cleveref}
\usepackage{thmtools}
\usepackage{thm-restate}

\numberwithin{equation}{section}

\newcommand{\E}{\mathbf{E}}

\DeclareMathOperator{\tr}{tr}

\newcommand{\T}{\mathsf{T}}

\newcommand{\floor}[1]{\lfloor #1 \rfloor}

\newcommand{\R}{\mathbb{R}}

\newcommand{\iid}{iid}

\newcommand{\opnorm}[1]{\| #1 \|_{\mathsf{op}}}

\newcommand{\hollow}[1]{\overset{\circ}{#1}}

\newtheorem{theorem}{Theorem} % same for example numbers
\newtheorem*{theorem*}{Theorem} % same for example numbers
\newtheorem{lemma}{Lemma}

\date{}

\title{\begin{center}
   An Elementary Proof of The Hanson-Wright Inequality
\end{center}}

\begin{document}

\maketitle

\vspace{-20pt}
\begin{abstract} 
The Hanson-Wright inequality establishes exponential concentration for quadratic forms $X^\T M X$, where $X$ is a vector with independent sub-Gaussian entries and with parameters depending on the Frobenius and operator norms of $M$. The most elementary proof to date is due to \citet{rudelson2013hanson}, who still rely on a convex decoupling argument due to \citet{bourgain1996random},  followed by Gaussian comparison to arrive at the result. In this note we sidestep this decoupling and provide an arguably simpler proof reliant only on elementary properties of sub-Gaussian variables and Gaussian rotational invariance. As a consequence we also obtain improved constants.
%that does the comparison directly by leveraging the martingale structure of the off-diagonal components of $X^\T M X$. We conclude by decoupling via rotational invariance.
\end{abstract}

\section{The Hanson-Wright Inequality}
Let $X_{1:n}$ be a sequence of mean zero, \iid-$\sigma^2$-sub-Gaussian random variables; $\E \exp \left(\lambda X_i \right) \leq \exp\left(\frac{\lambda^2 \sigma^2}{2}\right), \forall \lambda \in \R, i \in [n]$. In this note we prove the following exponential inequality.

\begin{theorem}
For every $\lambda \in  \left[0, \frac{1}{3c_2\opnorm{M}\sigma^2} \right) $  we have that:
    \begin{equation}
        \E \exp \left( \lambda \left[X^\T MX - \E [X^\T MX ]\right] \right)\leq \exp \left(c_1\lambda^2\sigma^4 \|M\|_F^2 \right).
    \end{equation}
    where we take $c_1 = 2$, $c_2=1$ if $M$ is diagonal-free and $c_1=20$, $c_2=4$ otherwise.
    
    Consequently for $t \geq 0$:
    \begin{equation}
        \mathbf{P}\left( \left|X^\T MX - \E [X^\T MX ]\right| \geq  t \right) \leq 2 \exp \left( -\left( \frac{t^2}{4c_1\sigma^4\|M\|_F^2} \wedge \frac{t}{6c_2\sigma^2 \opnorm{M}}\right)\right).
    \end{equation}
\end{theorem}

We let $M = (m_{ij})$ and define $\frac{M+M^\T}{2} \triangleq A  = (a_{ij})$. Observe that for any quadratic form $ x^\T M x= x^\T A x$ identically. Consequently, we have that
\begin{equation}
    X^\T M X = X^\T A X = \sum_{i=1}^n X_i^2 a_{ii} + 2 \sum_{j=2}^n X_j \sum_{i<j} a_{ij} X_i.
\end{equation}
The first term above is easy to analyze, since its just a sum of independent sub-exponential random variables. The second term is a little more tricky, and in the literature a convex decoupling inequality is typically used \citep{rudelson2013hanson}. Before we proceed, let us introduce $\hollow{A}$, the hollow of $A$, which is just a copy of $A$ but with its diagonal elements set to zero. Thus:
\begin{equation}\label{eq:hollow}
    X^\T \hollow A X = 2 \sum_{j=2}^n X_j \sum_{i<j} a_{ij} X_i.
\end{equation}
A quadratic form in the hollow of a symmetric matrix, such as $\sum_{j=2}^n X_j \sum_{i<j} a_{ij} X_i$, has a natural martingale structure that allows to directly produce a comparison inequality via repeated application of the tower rule. It is this "trick" that we refer to as elementary.

An approach similar in spirit to ours is due Latała \citep[see the appendix of][]{barthe2013transference} in which a decoupling inequality for U-statistics is used.  This idea is not dissimilar to decoupling \eqref{eq:hollow} using its martingale structure. We proceed to provide details of our direct approach below. 

\section{The Proof}

Let us introduce an auxiliary sequence $G_{1:n}$ of \iid\ Gaussian random variables with mean zero and variance $\sigma^2$. We have with $\E_n[\cdot] \triangleq \E[\cdot | X_{1:n-1}] $:
\begin{equation}
    \begin{aligned}
        &\E \exp \left(2\lambda \sum_{j=2}^n X_j \sum_{i<j} a_{ij} X_i  \right)\\
        & = \E \left[\exp \left(2\lambda \sum_{j=2}^{n-1} X_j \sum_{i<j} a_{ij} X_i \right) \E_n \exp \left(   2\lambda  X_n \sum_{i=2}^{n-1} a_{in} X_i \right)\right]\\
        &
        \leq  \E \left[\exp \left(2\lambda \sum_{j=2}^{n-1} X_j \sum_{i<j} a_{ij} X_i +   \frac{\sigma^2}{2} (2\lambda)^2  \left(\sum_{i=2}^{n-1} a_{in} X_i\right)^2 \right)\right] && (\textnormal{subG})\\
        &=  \E \left[\exp \left(2\lambda \sum_{j=2}^{n-1} X_j \sum_{i<j} a_{ij} X_i +  2\lambda G_n  \sum_{i=2}^{n-1} a_{in} X_i \right)\right] = (\dagger).
    \end{aligned}
\end{equation}
We can proceed similarly:
\begin{equation}
    \begin{aligned}
       (\dagger) &=\E \left[\exp \left(2\lambda \sum_{j=2}^{n-2} X_j \sum_{i<j} a_{ij} X_i +  2\lambda G_n  \sum_{i=2}^{n-2} a_{in} X_i \right)+2\lambda X_{n-1}(a_{n-1,n}G_n +\sum_{i=2}^{n-2} a_{i,n-1} X_i)  \right]\\
        &\leq \E \left[\exp \left(2\lambda \sum_{j=2}^{n-2} X_j \sum_{i<j} a_{ij} X_i +  2\lambda G_n  \sum_{i=2}^{n-2} a_{in} X_i \right)+\frac{(2\lambda)^2 \sigma^2 }{2}\left(a_{n-1,n}G_n +\sum_{i=2}^{n-2} a_{i,n-1} X_i\right)^2  \right]\\
        &=
        \E \left[\exp \left(2\lambda \sum_{j=2}^{n-2} X_j \sum_{i<j} a_{ij} X_i +  2\lambda G_n  \sum_{i=2}^{n-2} a_{in} X_i \right)+2\lambda G_{n-1} (a_{n-1,n}G_n +\sum_{i=2}^{n-2} a_{i,n-1} X_i)  \right]\\
        &\leq \dots = \E \exp \left(2\lambda \sum_{j=2}^n G_j \sum_{i<j} a_{ij} G_i  \right) = \E \exp \left( G^\T \hollow{A} G\right) = (\ddagger).
    \end{aligned}
\end{equation}
Indeed, the step $(\dots)$ can be established by combining $(\dagger)$, a finite induction argument and the following calculation:
\begin{equation}
    \begin{aligned}
        &\E \exp \left(\lambda  \begin{bmatrix}
            X_1 \\ G_{2:n}
        \end{bmatrix}^\T \begin{bmatrix}
            0 & a_{1,2:n}^\T\\
            a_{1,2:n} & \hollow{A}_{2:n,2:n}
        \end{bmatrix} \begin{bmatrix}
            X_1 \\ G_{2:n}
        \end{bmatrix} \right) \\
        & = \exp\left(2\lambda X_1 a_{1,2:n}^\T G_{2:n} + \lambda G_{2:n}^\T\hollow{A}_{2:n,2:n} G_{2:n} \right)\\
        &\leq \exp\left(\frac{\sigma^2}{2}(2\lambda)^2 (a_{1,2:n}^\T G_{2:n})^2 + \lambda G_{2:n}^\T\hollow{A}_{2:n,2:n} G_{2:n} \right)&& (\textnormal{subG})\\
        & = \exp\left(2\lambda G_1 a_{1,2:n}^\T G_{2:n} + \lambda G_{2:n}^\T\hollow{A}_{2:n,2:n} G_{2:n} \right) \\
        &=\E \exp \left(\lambda  \begin{bmatrix}
            G_1 \\ G_{2:n}
        \end{bmatrix}^\T \begin{bmatrix}
            0 & a_{1,2:n}^\T\\
            a_{1,2:n} & \hollow{A}_{2:n,2:n}
        \end{bmatrix} \begin{bmatrix}
            G_1 \\ G_{2:n}
        \end{bmatrix} \right).
    \end{aligned}
\end{equation}

Having established $(\ddagger)$, since $\hollow{A}$ is symmetric, we can write $G^\T \hollow{A}G = \sum_{i=1}^n \mu_i Z_i^2$ in distribution, where $\sigma Z_{1:n}$ is equal to $G_{1:n}$ in distribution relying on Gaussian rotational invariance. The next lemma is standard and bounds the moment generating function of this object.

\begin{lemma}\label{chi2lemma}
    Let $Z_{1:n}\sim N(0, \sigma^2 I_n)$. For every $\lambda \in \left[0, \frac{1}{3\max_{i\in[n]}|\mu_i|}\right]$ we have that:
    \begin{equation}
        \E \exp \left( \lambda \sum_{i=1}^{n}\mu_i Z_i^2 \right) \leq \exp\left(  \sum_{i=1}^{n} \lambda \mu_i +2\lambda^2 \mu_i^2 \right).
    \end{equation}
\end{lemma}
Note that $\max_{i\in[n]}|\mu_i|=  \sigma^2 \Big\|{\hollow{A}}\Big\|_{\mathsf{op}}$ and $\sum_{i=1}\mu_i^2 =\sigma^4  \Big\|\hollow{A} \Big\|_F^2$ in our case. Moreover, since $\hollow{A}$ is diagonal free $\lambda \sum_{i=1}^n \mu_i =\lambda \sigma^2 \tr \hollow{A} =0$. Hence we have the bound:
\begin{equation}
\label{eq:offdiag}
    (\ddagger) \leq \exp \left(2\lambda^2 \sigma^4 \Big\|{\hollow{A}}\Big\|_F^2 \right) \textnormal{ for }\lambda \in \left[0,\left(3\sigma^2\Big\|{\hollow{A}}\Big\|_{\mathsf{op}}\right)^{-1}\right].
\end{equation}
To analyze the diagonal terms we will require the following lemma.
\begin{lemma}\label{lem:bernstein}
    Let $X$ be $\sigma^2$-sub-Gaussian. We have that 
    \begin{equation}
       \E \exp \left( \lambda (X^2-\E X^2) \right) \leq \exp \left( 10 \lambda^2 \sigma^4 \right)
    \end{equation}
    for every nonnegative $\lambda$ satisfying $\lambda \leq \frac{1}{4\sigma^2 }$ .
\end{lemma}
We proceed to apply the above lemma.  On the region $ \{\lambda : \max |4\lambda a_{ii}|\sigma^2 <1 \}$ we have that 
\begin{equation}
\label{eq:diag}
    \begin{aligned}
        \E \exp \left( \lambda \sum_{i=1}^n (X_i^2-\E Xi^2) a_{ii}\right) &= \prod_{i=1}^n \E \exp \left( \lambda  (X_i^2-\E Xi^2) a_{ii}\right)\\
        & \leq  \prod_{i=1}^n  \exp \left(   10 \lambda^2 \sigma^4 a_{ii}^2\right)\\
        &\leq \exp \left(   10 \lambda^2 \sigma^4 \left\| A-\hollow A\right\|_F^2\right).
    \end{aligned}
\end{equation}
To finish the proof, we combine \eqref{eq:offdiag} and \eqref{eq:diag} with the Cauchy-Schwarz inequality (noting that this is unnecessary if $M$ is diagonal-free):
\begin{multline}
    \E \exp \left( \lambda \left[ X^\T M X - \E [X^\T MX]\right] \right)
    =\E \exp \left( \lambda \left[ X^\T A X - \E [X^\T AX]\right] \right)
    \\
    \leq \sqrt{\E \exp \left( 2\lambda \left[ X^T (A-\hollow{A}) X - \E [X^\T (A-\hollow {A})X]\right] \right)\E \exp \left( 2\lambda \left[ X^\T \hollow{A} X - \E [X^\T \hollow{A}X]\right] \right)}
    \\
    \leq \sqrt{  \exp \left(   40\lambda^2 \sigma^4 \left\| A-\hollow A\right\|_F^2\right)  \exp \left(8\lambda^2 \sigma^4 \Big\|{\hollow{A}}\Big\|_F^2 \right)  }
    \\
    \leq  \exp \left(   20\lambda^2 \sigma^4 \left\| A-\hollow A\right\|_F^2+ 20\lambda^2 \sigma^4 \Big\|{\hollow{A}}\Big\|_F^2 \right) 
\end{multline}
as long as  $\max |4(2\lambda) a_{ii}|\sigma^2 <1$ and $\left|3(2\lambda)\left\| \hollow{A}\right\|_{\mathsf{op}}\right|\sigma^2 <1$. The result follows since $\left\| \hollow{A}\right\|_{\mathsf{op}} \leq 2 \left\| A\right\|_{\mathsf{op}}.$

\section{Proofs of Auxiliary Lemmata}

\begin{proof}[Proof of \Cref{chi2lemma}]
For $\lambda \leq \frac{1}{3\max_{i\in[n]}  |\mu_i|}$ we have that:
    \begin{equation}
        \begin{aligned}
             \E \exp \left( \lambda \sum_{i=1}^{n}\mu_i Z_i^2 \right)&=\E \exp \left(- \frac{1}{2}  \sum_{i=1}^{n} \log (1- 2\lambda \mu_i ) \right) \\
             &\leq \E \exp \left( \sum_{i=1}^{n} \lambda \mu_i   +2\lambda^2 \mu_i^2\right) &&(-\log(1-2x) \leq 2x+4x^2, |x| \leq 1/3).
        \end{aligned}
    \end{equation}
\end{proof}

\begin{proof}[Proof of \Cref{lem:bernstein}]
We proceed by expanding the moment generating function.
    \begin{equation}
        \begin{aligned}
            \E \exp \left( \lambda (X^2-\E X^2) \right)
            &= 1+\sum_{k=2}^\infty \frac{ \lambda^k \E (X^2-\E X^2)^k}{k!}\\
            %!
            & \leq  1+\sum_{k=2}^\infty  (2\lambda\sigma^{2})^k \cosh\left(\frac{1}{2}\right) &&(\dagger)\\
            &= 1+ \cosh\left(\frac{1}{2}\right)(2\lambda\sigma^{2})^2\sum_{k=0}^\infty  (2\lambda\sigma^{2})^k  \\
            &
            \leq 1 + \cosh\left(\frac{1}{2}\right) \frac{(2\lambda\sigma^{2})^2}{1-2\lambda\sigma^2}
        \end{aligned}
    \end{equation}
valid on the region $ \{\lambda : |2\lambda \sigma^2 | <1 \}$. In particular for, nonnegative $\lambda \leq \frac{1}{4\sigma^2 }$ we have that
\begin{equation}
    \begin{aligned}
          \E \exp \left( \lambda (X^2-\E X^2) \right) \leq  1+10\lambda^2 \sigma^4\leq \exp \left( 10 \lambda^2 \sigma^4 \right)
    \end{aligned}
\end{equation}
as was required. The step $(\dagger)$ can be shown as follows:

\begin{equation}
    \begin{aligned}
        \E (X^2-\E X^2)^k &= \sum_{j=0}^k {k \choose j } (-1)^{k-j} \left( \E X^2\right)^{k-j} \E X^{2j}&&(\textnormal{Binomial Theorem})\\
        &\leq \sum_{j=0, \: j-k \:\mathrm{even} }^k  {k \choose j } \left( \E X^2\right)^{k-j} \E X^{2j} 
        \\
        &\leq \sum_{j=0, \: j-k \:\mathrm{even} }^k  {k \choose j } \left( \E \sigma^2\right)^{k-j}  \left(2^j j! \sigma^{2j} \right) &&(\textnormal{sub-Gaussian moments})
        \\
        &= \sigma^{2k} \sum_{j=0, \: j-k \:\mathrm{even} }^k  {k \choose j }  \left(2^j j! \right)  
        \\
        &
        =\sigma^{2k}\sum_{m=0 }^{\floor {k}}  {k \choose k-2m }  2^{k-2m} (k-2m)!  && (j=k-2m)
        \\
         &
        =(2\sigma^{2})^kk! \sum_{m=0 }^{\floor {k}} \frac{1}{(2m)!} \left(\frac{1}{2}\right)^{2m}
        \\
         &\leq (2\sigma^{2})^kk! \sum_{m=0 }^{\infty} \frac{1}{(2m)!} \left(\frac{1}{2}\right)^{2m}
        \\
        &
        =  \cosh\left(\frac{1}{2}\right) (2\sigma^{2})^kk!
    \end{aligned}
\end{equation}
which finishes the proof.
\end{proof}

%\section*{Acknowledgments}

\addcontentsline{toc}{section}{References}

\bibliographystyle{tmlr}
\bibliography{main.bib}

\end{document}